\documentstyle[11pt,twoside]{article}

\newcommand{\dopu}{{:}\allowbreak\ }
\newcommand{\Id}{I\mkern-1mud}
\newcommand{\eps}{\varepsilon}
\newcommand{\del}{\delta}
\newcommand{\settminus}{{\setminus}}
\newcommand{\sst}{\scriptstyle}
\newcommand{\maxp}{\mathop{\rm \vphantom{p}max}} 

\newcommand{\N}{{\bf N}}

\newcommand{\proof}{\par\noindent{\em Proof. }}
\newcommand{\eop}{\nopagebreak\hspace*{\fill}$\Box$}

\newtheorem{theo}{Theorem}
\newtheorem{lemma}[theo]{Lemma}
\newtheorem{cor}[theo]{Corollary}
\newtheorem{prop}[theo]{Proposition}

\newcommand{\bea}{\begin{eqnarray*}}
\newcommand{\eea}{\end{eqnarray*}}

\newcommand{\muu}[1]{\mu_{#1}(\{#1\})}
\newcommand{\muuu}[2]{\mu_{#1}(\{#2\})}

\newif\ifrefsc
\let\thebibliographyalt=\thebibliography                                %
\def\thebibliography#1                                                  %
 {\ifrefsc
 \section*{\normalsize\hfil\sc References\hfill}
 \small
 \list                    %
 {[\arabic{enumi}]}{\settowidth\labelwidth{[#1]}\leftmargin\labelwidth  %
 \advance\leftmargin\labelsep                                           %
 \usecounter{enumi}}                                                    %
 \def\newblock{\hskip .11em plus .33em minus .07em}                     %
 \sloppy\clubpenalty4000\widowpenalty4000                               %
 \sfcode`\.=1000\relax                                                  %
 \else\thebibliographyalt{#1}\fi}                                         %

\refsctrue

\begin{document}

\begin{center}
{\Large\bf 
An elementary approach to the \\
Daugavet equation
}
\\[22 pt]
{\sc  Dirk Werner}
\end{center}
\bigskip

\begin{quote}
\small
{\sc Abstract.}		
Let $T\dopu C(S)\to C(S)$ be a bounded linear operator. We present a
necessary and sufficient condition for the so-called Daugavet
equation
$$ \|\Id+T\| = 1+\|T\| $$
to hold, and we apply it to weakly compact operators and to
operators factoring through $c_{0}$. Thus we obtain very simple proofs of
results by Foia\c{s}, Singer, Pe\l czy\'nski, Holub and others.
\end{quote}
\bigskip

\noindent
If $E$ is a real Banach space, let us say that an operator 
$T\dopu E\to E$ satisfies the {\em Daugavet equation\/}
if
$$
\|\Id+T\| = 1+\|T\|.
$$
Daugavet 
\cite{Daug} 
proved that every compact operator $T\dopu C[0,1]\to
C[0,1]$ satisfies this equation, and Foia\c{s} and Singer 
\cite{FoiSin}
extended his result to weakly compact operators. Later, these
theorems were rediscovered by Kamowitz 
\cite{Kamow} 
and Holub 
\cite{Holub1}. 
Pe\l{}czy\'nski
observed that the Foia\c{s}-Singer argument can be used to prove the
Daugavet equation for weakly compact 
operators on a $C(S)$-space provided $S$ has no isolated points,
cf.\ 
\cite[p.~446]{FoiSin}. 
(This restriction on $S$ is easily seen to be necessary.)
Holub also showed, with no assumption on $T$, that $T$ or $-T$
fulfills the Daugavet equation for which Abramovich 
\cite{Abra1} 
gave another proof
valid for general $S$ rather than the unit interval. Yet another
argument was suggested in 
\cite[p.~343]{HWW}.
Dually, a number of authors have investigated the Daugavet equation
for operators on an $L^1$-space; for more precise references we refer
to 
\cite{Abra2} 
and 
\cite{Schm}. 
In the latter paper Schmidt proved the Daugavet equation
for weakly compact operators on an atomless $L^1$-space using Banach
lattice techniques.
(Actually, his result is more general than that.)
A different class of operators was considered by Holub
\cite{Holub2}
who showed the  Daugavet equation for the ideal of operators on
$C[0,1]$ factoring through $c_{0}$; Ansari 
\cite{Ansari} 
generalised this result to such operators on $C(S)$-spaces, where $S$
has no isolated points.

In this paper we suggest a unified and elementary approach to all the
results just mentioned. Our basic idea is to represent an operator
$T\dopu C(S)\to C(S)$ by its stochastic kernel, that is the family 
of measures $(\mu_s)_{s\in S}$ defined by $\mu_s=T^*\del_s$; i.e.,
$$
\int_S	f\,d\mu_{s} = \langle f,\mu_s \rangle =\langle Tf,\del_s \rangle
= (Tf)(s).
$$
We then have $\|T\|=\sup_s \|\mu_s\|$, and the function $s\mapsto
\mu_s$ is continuous for the weak$^*$ topology of $M(S)\cong
C(S)^*$. The operator $T$ is weakly compact if and only if $s\mapsto
\mu_s$ is continuous for the weak topology of $M(S)$ (meaning the
$\sigma(M(S),M(S)^*)$-topology), and $T$ is compact if and only if $s\mapsto
\mu_s$ is norm continuous; see 
\cite[p.~490]{DS}.
Note that the identity operator is represented by the family of Dirac
measures $(\del_s)_{s\in S}$.

A different approach to the Daugavet equation for weakly compact
operators was taken by
Abramovich, Aliprantis, and Burkinshaw 
\cite{AbraAB} 
who used ideas from Banach lattice theory, and  Ansari 
\cite{Ansari}
was able to incorporate Holub's result on
$c_{0}$-factorable operators into their scheme. However, these
arguments seem to be less elementary than the very simple
calculations presented here. 

We finally mention the recent papers 
\cite{Abra2}
and 
\cite{Woj92} 
whose results are not covered by this note.
\smallskip

\noindent
{\em Acknowledgement.} I would like to thank the referee for pointing
out Ansari's paper to me.

\begin{prop}\label{one}
If $S$ is a compact Hausdorff space and $T\dopu C(S)\to C(S)$ is a
bounded linear operator, then
$$
\max\{\|\Id + T\|, \|\Id-T\| \} = 1+\|T\|.
$$
\end{prop}

\proof
Let $(\mu_s)_{s\in S}$ be the representing kernel of $T$. Then
\bea
\maxp_\pm \|\Id \pm T\| &=&
\maxp_\pm \sup_{s\in S} \|\del_s \pm \mu_s\|  \\
&=&
\sup_{s\in S} \maxp_\pm  \bigl( |\del_s \pm \mu_s|(\{s\}) +
|\del_s \pm \mu_s|(S\settminus \{s\}) \bigr) \\
&=&
\sup_{s\in S} \maxp_\pm \bigl( |1\pm \muu{s}| + |\mu_s|(S\settminus \{s\})
\bigr) \\
&=&
\sup_{s\in S}  \bigl( 1 + |\muu{s}| + |\mu_s|(S\settminus \{s\})
\bigr) \\
&=&
\sup_{s\in S} \bigl( 1 + \|\mu_s\| \bigr)
~=~ 1+\|T\|.
\eea
\eop

\begin{cor}\label{two}
If $E$ is an (AL)-space or an (AM)-space and $T\dopu E\to E$ is a
bounded linear operator, then
$$
\max\{\|\Id + T\|, \|\Id-T\| \} = 1+\|T\|.
$$
\end{cor}

\proof
An (AL)-space $E$ is representable as $L^1(\mu)$ for some localisable
measure $\mu$, hence $E^*$ is representable as $L^{\infty}(\mu) \cong
C(S)$. So the assertion 
follows from Proposition~\ref{one} by passing to $T^*$. If $E$ is an
(AM)-space, then $E^*$ is an (AL)-space, and again we obtain the assertion
by considering the adjoint operator.
\eop
\bigskip

We now formulate a technical condition that will allow us to prove the
Daugavet equation for weakly compact operators and for
$c_{0}$-factorable operators.

\begin{lemma}\label{5}
Let $S$ be a compact Hausdorff space and $T\dopu C(S)\to C(S)$  a
bounded linear operator with representing kernel $(\mu_{s})_{s\in
S}$. If the kernel satisfies
$$
\sup\limits_{s\in U} \muu{s} \ge0 \mbox{ for all nonvoid
open sets }U\subset S, \eqno(*)
$$
then
$$
\|\Id + T\| = 1+\|T\|.
$$
In fact, a necessary and sufficient condition for this to hold is
$$
\sup_{\{s\dopu \|\mu_{s}\| > \|T\|-\eps \}} \muu{s} \ge 0
\qquad\forall \eps>0  .   \eqno(**)
$$
\end{lemma}

\proof
We have
$$
\|\Id+T\| = \sup_{s\in S} \|\del_s+\mu_s\| =
\sup_{s\in S} \bigl( |1+\muu{s}| + |\mu_s|(S\settminus\{s\}) \bigr)
$$
and
$$
1+\|T\| = \sup_{s\in S} \bigl( 1+ \|\mu_s\| \bigr)
=\sup_{s\in S} \bigl( 1+ |\muu{s}| + |\mu_s|(S\settminus\{s\}) \bigr);
$$
so problems with showing the Daugavet equation can only arise in case 
some of the $\muu{s}$ are negative. 

Given $\eps>0$, we now apply 
$(*)$ to the open set
$U=\{s\in S\dopu \|\mu_s\|>\|T\|-\eps\}$ 
(that is, we apply $(**)$)
and obtain
\bea
\|\Id+T\| &\ge&
\sup_{s\in U} \|\del_s + \mu_s\| \\
&=&
\sup_{s\in U} \bigl( |1+ \muu{s}| + |\mu_s|(S\settminus\{s\}) \bigr) \\
&\ge&
\sup_{\sst s\in U \atop \sst \muu{s}\ge -\eps}
\bigl( 1+\|\mu_s\| + \muu{s} -|\muu{s}| \bigr) \\
&\ge&
1+\|T\|-\eps +
\sup_{\sst s\in U \atop \sst \muu{s}\ge -\eps}
\bigl( \muu{s} -|\muu{s}| \bigr) \\
&\ge&
1+\|T\| -3\eps;
\eea
hence $T$ satisfies the Daugavet equation.

A similar calculation shows that $(**)$ is not only sufficient, but
also necessary.
\eop
\bigskip

Next, we deal with weakly compact operators.

\begin{lemma}\label{6}
If $S$ is a compact Hausdorff space without isolated points and
$T\dopu C(S)\to C(S)$ is weakly compact, then $T$ fulfills $(*)$ of
Lemma~\ref{5}.
\end{lemma}

\proof
To prove this lemma we argue by contradiction. Suppose there is a
nonvoid open set $U\subset S$ and some $\beta>0$ such that
$$
\muu{s} <-2\beta \qquad \forall s\in U.
$$
At this stage we note that, for each $t\in S$, the function 
$s\mapsto \muuu{s}{t}$
is continuous, since $T$ is weakly compact. For $\mu\mapsto \mu(\{t\})$
is in $M(S)^*$ and, as noted in the introduction, $s\mapsto \mu_s$
is weakly continuous.

Returning to our argument we pick some $s_0\in U$ and consider the set
$$
U_1 = \{s\in U\dopu |\muuu{s}{s_0} - \muu{s_0}|<\beta \}
$$
which---as we have just observed---is an open neighbourhood of $s_0$.
Since $s_0$ is not isolated, there is some $s_1\in U_1$, $s_1\neq s_0$.
We thus have
$$
\muu{s_1} <-2\beta,
$$
because $s_1\in U$, and
$$
\muuu{s_1}{s_0} < \muu{s_0} + \beta < -2\beta + \beta = -\beta.
$$
In the next step we let
$$
U_2 = \{s\in U_1\dopu |\muuu{s}{s_1} - \muu{s_1}|<\beta \} \ (\subset U).
$$
Likewise, this is an open neighbourhood of $s_1$, hence there is some
$s_2\in U_2$, $s_2\neq s_1$, $s_2\neq s_0$. We conclude, using
that $s_2\in U$, $s_2\in U_2$ and $s_2\in U_1$,
\bea
\muuu{s_2}{s_2}  &<&   -2\beta , \\
\muuu{s_2}{s_1}  &<&   -\beta , \\
\muuu{s_2}{s_0}  &<&   -\beta .
\eea
Thus we inductively define a descending sequence of 
open sets $U_n\subset U$ and distinct points
$s_n\in U$ by
\bea
U_{n+1} &=& \{s\in U_n\dopu |\muuu{s}{s_n} - \muu{s_n}| <\beta\}, \\
s_{n+1} &\in& U_{n+1}\settminus \{s_0,\ldots,s_n\}
\eea
yielding
$$
\muuu{s_n}{s_j} < -\beta \qquad \forall j=0,\ldots,n-1.
$$
Consequently,
$$
\|T\|\ge \|\mu_{s_n}\| \ge |\mu_{s_n}|(\{s_0,\ldots,s_{n-1}\}) \ge n\beta
\qquad\forall n\in \N,
$$
which furnishes a contradiction.
\eop
\bigskip

Lemmas~\ref{5} and \ref{6} immediately yield the first main result of this
note.

\begin{theo}\label{three}
Suppose $S$ is a compact Hausdorff space without isolated points.
If $T\dopu C(S)\to C(S)$ is weakly compact, then
$$
\|\Id + T\|  = 1+\|T\|.
$$
\end{theo}

%

\begin{cor}\label{four}
If $\mu$ is an atomless measure and $T\dopu L^1(\mu)\to L^1(\mu)$
is weakly compact, then
$$
\|Id+T\| = 1 + \|T\|.
$$
\end{cor}

\proof
By changing measures if necessary we may assume that $L^1(\mu)^*\cong
L^\infty(\mu)$ canonically. (If $L^1(\mu)\cong L^1(\nu)$ and $\mu$ is
atomless, then so is $\nu$, since atomless measure spaces are 
characterised by the fact that the unit balls of the corresponding
$L^1$-spaces fail to possess extreme points.) Now $L^\infty(\mu)$
is isometric to some $C(S)$-space, where $S$ does not contain any 
isolated point. It remains to observe that $T^*$ is weakly compact as well
\cite[p.~485]{DS} 
and to apply Theorem~\ref{three}.
\eop
\bigskip

\noindent
{\em Remarks.}
(1)
If $T$ is compact, the proof of Lemma~\ref{6} can considerably be 
simplified. In fact, if $\muu{s}<-2\beta <0$ on an open nonvoid set $U$,
let us pick some $s\in U$ and consider the set 
$$
U_1 = \{t\in U\dopu \|\mu_s - \mu_t\|<\beta\}.
$$
Since $T$ is compact, this is an open neighbourhood of $s$, and for
each $t\in U_1$ we deduce that
$$
\muuu{s}{t} \le
\muu{t} + |\muu{t} - \muuu{s}{t}|  
< -2\beta + \|\mu_t-\mu_s\|
< -\beta.
$$
Since $s$ is not isolated, there are infinitely many distinct points
$t_1, t_2, \ldots \in U_1$, and we obtain
$|\mu_s|(\{t_1,t_2,\ldots\}) =\infty$,
a contradiction.

(2)
The proof of Theorem~\ref{three} shows that weakly compact operators
on $C_0(S)$, $S$ locally compact without isolated points, satisfy
the Daugavet equation.

(3)
We also see immediately that positive operators on $C(S)$-spaces
(and likewise on (AL)- and (AM)-spaces) satisfy the Daugavet equation.

(4)
For weakly compact operators $T$ on $C(S)$, represented by 
$(\mu_s)_{s\in S}$, the functions $\varphi_A \dopu
s\mapsto \mu_s(A)$, $A\subset S$
a Borel set, are continuous; in fact, weakly compact operators 
are characterised by this property 
\cite[p.~493]{DS}. 
In Lemma~\ref{6} it is even enough to assume that
only the functions $\varphi_{\{t\}}$, $t\in S$, are continuous,
provided $S$ has no isolated points.
Hence also such operators satisfy the Daugavet equation.
A special case of this situation (a trivial one, though) occurs
if $\mu_s(\{t\})=0$ for all $s,t\in S$; see also the following remark.

(5)
A particular class of operators for which $(*)$ of Lemma~\ref{5} is
valid are those for which 
$$
\{t\in S\dopu \muuu{s}{t}=0\ \forall s\in S\} \mbox{ is dense in $S$.}
\eqno(*{*}*)
$$ 
Since this class is seen to contain
the {\em almost diffuse operators\/}
of Foia\c{s} and Singer,
we have obtained their result that almost
diffuse operators satisfy the Daugavet equation.

\bigskip
This last remark easily leads to Ansari's extension of Holub's theorem that
operators on $C[0,1]$ factoring through $c_{0}$ satisfy the
Daugavet equation.

\begin{theo}\label{seven}
If $S$ is a compact Hausdorff space without isolated points and
$T\dopu C(S)\to C(S)$ factors through $c_{0}$, then 
$$
\|\Id+T\| = 1+\|T\|.
$$
\end{theo}

\proof        
Let $(\mu_s)_{s\in S}$ be the representing kernel of $T$. 
By remark~(5) it is enough to show that
$$
S' := \{t\in S\dopu \muuu{s}{t}=0 \ \forall s\in S\}
$$ 
is dense in $S$. 
Let us write $T=T_{2}T_{1}$ with bounded linear operators $T_{1}\dopu
C(S)\to c_{0}$, $T_{2}\dopu c_{0}\to C(S)$. We have
\bea
(T_{1}f)(n) 
&=& 
\int_{S} f\, d\rho_{n} \qquad\forall n\in\N,  \\
\bigl( T_{2}(a_{n}) \bigr)(s) 
&=& 
\sum_{n=1}^{\infty} \nu_{s}(n)a_{n}
\qquad\forall s\in S
\eea
for a sequence of measures $\rho_{n}$ and a family $\bigl( \nu_{s}(n)
\bigr)_{n}$ of sequences in $\ell_{1}$. Consequently,
$$
\mu_{s} = \sum_{n=1} ^{\infty} \nu_{s}(n) \rho_{n}.
$$
Now $S' \supset \bigcap_{n} \{t\in S\dopu \rho_n(\{t\})=0\}$, which
is a set whose complement is at most countable. Since no point in $S$
is isolated, countable sets are of the first category, and Baire's
theorem implies that $S'$ is dense.
\eop

\bigskip
\noindent
{\em More remarks.}
(6)
The same proof applies to operators that factor through a
$C(K)$-space where $K$ is a countable 
compact space, since on such spaces all regular Borel measures
are discrete. We recall that there are countable compact spaces $K$
such that $C(K)$ is not isomorphic to $c_{0}$.

(7)
The Baire argument in Theorem~\ref{seven} implies a very simple proof
of Theorem~\ref{three} if in addition $S$ is supposed to be
separable. In fact, let us show that then $(*{*}*)$ of Remark~(5) holds.
The complement of the set spelt out there is $\{t\in S\dopu \exists
s\in S\ \muuu{s}{t}\neq0 \}$. Since $s\mapsto \muuu{s}{t}$ is
continuous, this is, with $\{s_{1},s_{2},\ldots\}$ denoting a
countable dense subset of $S$, $\bigcup_{n} \{t\in S\dopu 
\muuu{s_{n}}{t}\neq0 \}$ and hence a countable union of countable
sets, i.e., of the first category. Again,
$\{t\in S\dopu  \muuu{s}{t}=0\ \forall s\in S \}$ must be dense.

(8)
We finally wish to comment on the case of complex scalars. All the
results and proofs in this paper remain valid---mutatis mutandis---in
the setting of complex Banach spaces. In Proposition~\ref{one} the
proper formulation of the conclusion is 
$$
\max_{|\lambda|=1} \|\Id+\lambda T\| = 1 + \|T\|,
$$
and $(*)$ in Lemma~\ref{5} should be replaced by
$$
\sup_{s\in U} \bigl( |1+\muu{s}| - (1+|\muu{s}|) \bigr) \ge 0
\mbox{ \em for all nonvoid  open sets }U\subset S.
$$


\bigskip
\bigskip
\noindent
I.\ Mathematisches Institut
\\
 Freie Universit\"{a}t Berlin
\\
Arnimallee 3
\\
D-14\,195 Berlin 
\\
Federal Republic of Germany
\\[1mm]
e-mail: werner@math.fu-berlin.de

\end{document}